\documentclass[12pt]{article}

\usepackage{amsmath}
\usepackage{theorem}
\usepackage{amssymb}
\usepackage{latexsym}
\usepackage{a4wide}
\usepackage{youngtab}

\newtheorem{thm}{Theorem}
\newtheorem{lem}[thm]{Lemma}
\newtheorem{prop}[thm]{Proposition}
\newtheorem{cor}[thm]{Corollary}

{\theorembodyfont{\rmfamily} \newtheorem{exa}[thm]{Example}}
\newenvironment{rem}{\noindent{\bf Remark.}}{\newline}
\newenvironment{pf}{\noindent{\bf Proof.}}{\hbox{}\hfill $\Box$}

\newcommand{\Q}{\mathbb{Q}}
\newcommand{\R}{\mathbb{R}}
\newcommand{\Z}{\mathbb{Z}}

\newcommand{\B}{\mathcal{B}}

\newcommand{\LL}{\mathcal{L}}

\newcommand{\LLlam}{\mathcal{L}(\lambda)}

\newcommand{\Plexg}{{>_{\rm lex}}}
\newcommand{\Plexs}{{<_{\rm lex}}}

\newcommand{\Ftilde}{\widetilde{F}}

\begin{document}

\title{An algorithm to compute the canonical basis of an irreducible 
$U_q(\mathfrak{g})$-module}
\author{Willem A. de Graaf,  University of Utrecht, The Netherlands\\
email: {\tt degraaf@math.uu.nl}}
\date{}
\maketitle

\begin{abstract}
We describe an algorithm to compute the canonical basis of an irreducible module
over a quantized enveloping algebra of a finite-dimensional semisimple Lie algebra.
The algorithm works for modules that are constructed as a submodule of a tensor product
of modules with known canonical bases.
\end{abstract}

\section{Introduction}

In this paper I consider the problem of constructing the canonical basis 
(cf. \cite{lusztig}) of 
an irreducible module over a quantized enveloping algebra. There are several possible
ways to approach this problem, and they may depend on how the module is constructed.
In \cite{gra8} an algorithm is described that works for any module, provided that
we have a method for computing the action of elements of the algebra. In \cite{lectof},
\cite{lecouvey} the irreducible module is first constructed as a submodule of a
tensor product of other modules. Then using the known canonical bases
of these other modules, an algorithm is described for constructing the canonical 
basis of the submodule.\par
Since constructing irreducible modules as submodules of tensor products can be
rather efficient (cf. \cite{gra9}), it is worthwhile to have an algorithm that is 
tailored to this situation. Therefore, 
in this paper I take the second approach above. In fact, I will describe an algorithm 
that is very similar to the ones in \cite{lectof}, \cite{lecouvey}. The main difference
is that I do not assume that the root system is of a certain type. The algorithm given
here works for all types, assuming that somehow we know the canonical bases of the
fundamental modules. These can, for instance, be constructed using the algorithm of 
\cite{gra8}.\par
This paper is organised as follows. In Section \ref{prelim} the theoretical concepts,
and notation that we use are introduced. Then in Section \ref{psi} a result
is described concerning the form of the elements of the canonical basis of a tensor 
product. Then in Section \ref{monomial} this is used, along with the description
of a monomial basis of an irreducible module (from \cite{lakshmibai}), to give an
algorithm for constructing the canonical basis. Finally in Section \ref{an} this
algorithm is compared to the algorithm from \cite{lectof} in the $A_n$-case. It is
shown that in that case both algorithms are very similar (but not exactly the same).

\section{Preliminaries}\label{prelim}

In this section we briefly sketch the concepts and notation that we will be using.
Our main reference is \cite{jcj}.\par
Let $\mathfrak{g}$ be a semisimple Lie algebra over $\mathbb{C}$. By $\Phi$ we
denote the root system of $\mathfrak{g}$, and $\Delta=\{\alpha_1,\ldots,\alpha_l\}$
will be a fixed set of simple roots of $\Phi$. Let $W$ denote the Weyl group of $\Phi$,
which is generated by the simple reflections $s_i=s_{\alpha_i}$ for $1\leq i\leq l$.
Let $\R \Phi$ be the vector space over $\R$ spanned by $\Phi$. On $\R\Phi$ we fix
a $W$-invariant inner product $(~,~)$ such that $(\alpha,\alpha)=2$ for short roots
$\alpha$. This means that $(\alpha,\alpha)=2,4,6$ for $\alpha\in\Phi$. \par
We work over the field $\Q(q)$. For $\alpha\in\Phi$ set $q_{\alpha}=q^{(\alpha,\alpha)/2}$.
For $n\in\Z$ we set $[n]_{\alpha} = q_{\alpha}^{-n+1} + q_{\alpha}^{-n+3}+\cdots +
q_{\alpha}^{n-1}$. Also $[n]_{\alpha}! = [n]_{\alpha}[n-1]_{\alpha}\cdots [1]_{\alpha}$ and
$$\begin{bmatrix} n \\ k \end{bmatrix}_{\alpha} = 
\frac{[n]_{\alpha}!}{[k]_{\alpha}![n-k]_{\alpha}!}.$$
Let $\Delta=\{\alpha_1,\ldots, \alpha_l\}$ be a simple system of $\Phi$.
Then the quantized enveloping algebra $U_q=U_q(\mathfrak{g})$ is the associative
algebra (with one) over $\Q(q)$ generated by $F_{\alpha}$, $K_{\alpha}$,
$K_{\alpha}^{-1}$, $E_{\alpha}$ for $\alpha\in\Delta$, subject to the following
relations
\begin{align*}
K_{\alpha}K_{\alpha}^{-1} &= K_{\alpha}^{-1}K_{\alpha} = 1,~
K_{\alpha}K_{\beta} = K_{\beta}K_{\alpha}\\
E_{\beta} K_{\alpha} &= q^{-(\alpha,\beta)}K_{\alpha} E_{\beta}\\
K_{\alpha} F_{\beta} &=  q^{-(\alpha,\beta)}F_{\beta}K_{\alpha}\\
E_{\alpha} F_{\beta} &= F_{\beta}E_{\alpha} +\delta_{\alpha,\beta}
\frac{K_{\alpha}-K_{\alpha}^{-1}}{q_{\alpha}-q_{\alpha}^{-1}} \\
\sum_{k=0}^{1-\langle \beta,\alpha^{\vee}\rangle } & (-1)^k 
\begin{bmatrix} 1-\langle \beta,\alpha^{\vee}\rangle \\ k 
\end{bmatrix}_{\alpha} E_{\alpha}^{1-\langle \beta,\alpha^{\vee}\rangle-k}
E_{\beta} E_{\alpha}^k = 0\\
\sum_{k=0}^{1-\langle \beta,\alpha^{\vee}\rangle } & (-1)^k 
\begin{bmatrix} 1-\langle \beta,\alpha^{\vee}\rangle \\ k 
\end{bmatrix}_{\alpha} F_{\alpha}^{1-\langle \beta,\alpha^{\vee}\rangle-k}
F_{\beta} F_{\alpha}^k= 0,
\end{align*}
where the last two relations are for all $\alpha\neq \beta$. \par
Let $U^-$, $U^0$, $U^+$ be the subalgebras of $U_q$ generated by respectively, 
$F_{\alpha}$ for $\alpha\in\Delta$, $K_{\alpha}^{\pm 1}$ for $\alpha\in\Delta$, and
$E_{\alpha}$ for $\alpha\in\Delta$. Then as a vector space $U_q\cong U^-\otimes
U^0 \otimes U^+$ (\cite{jcj}, Theorem 4.21). Let $\nu=\sum_k a_k \alpha_k$
with $a_k\in \Z_{\geq 0}$. Then we let $U^+_{\nu}$ ($U_{\nu}^-$) be the subspace of 
$U^+$ ($U^-$) spanned by all $E_{\alpha_{i_1}}\cdots E_{\alpha_{i_r}}$
($F_{\alpha_{i_1}}\cdots F_{\alpha_{i_r}}$) such that $\alpha_{i_1}
+\cdots +\alpha_{i_r}=\nu$. \par
We let $\lambda_1,\ldots, \lambda_l$ denote the fundamental weights, and 
$P=\Z \lambda_1+\cdots +\Z \lambda_l$ is the weight lattice. Also 
$P^+ = \Z_{\geq 0}\lambda_1+\cdots +\Z_{\geq 0}\lambda_l$ is the set of dominant 
weights. Now for every dominant $\lambda\in P^+$ there is an irreducible 
$U_q$-module $V(\lambda)$. We have that $V(\lambda)$ is spanned by vectors $v_{\mu}$
for $\mu\in P$, with $K_{\alpha} \cdot v_{\mu} = q^{(\mu,\alpha)} v_{\mu}$. These
$v_{\mu}$ are called weight-vectors of weight $\mu$. Among them there is the vector
$v_{\lambda}$ (which is unique upto scalar multiples), with $U^+\cdot v_{\lambda}=0$.
This $v_{\lambda}$ is called the highest-weight vector. We have that $V(\lambda) = 
U^-\cdot v_{\lambda}$. Furthermore, every finite-dimensional irreducible $U_q$-module 
is isomorphic to a $V(\lambda)$ (\cite{jcj}, Theorem 5.10). \par
Let $M$ be a finite-dimensional $U_q$-module. Then $M$ has a crystal base 
$(\mathcal{M}, \B)$ as defined in \cite{jcj}, 9.4. Here $\mathcal{M}$ is an
$A$-submodule of $M$, where $A$ is the subring of $\Q(q)$ consisting of rational
functions without pole at $0$. And $\B$ is a basis of $\mathcal{M}/q\mathcal{M}$. 
For $\alpha\in\Delta$, we have the Kashiwara operators
$\Ftilde_{\alpha}, \widetilde{E}_{\alpha} : \mathcal{M}\to \mathcal{M}$, and
the induced operators $\Ftilde_{\alpha}, \widetilde{E}_{\alpha} : \B\to\B\cup\{0\}$
(\cite{jcj}, 9.2, 9.4). \par
There is a $\Q$-algebra isomorphism $\overline{\phantom{a}} : U_q\to U_q$
with $\overline{q}=q^{-1}$, $\overline{E}_{\alpha}=E_{\alpha}$,
$\overline{F}_{\alpha}=F_{\alpha}$, and  $\overline{K}_{\alpha}=K^{-1}_{\alpha}$
(\cite{jcj}, Proposition 11.9). If $V(\lambda)$ is an irreducible $U_q$-module
with highest weight $\lambda$, and fixed highest-weight vector $v_{\lambda}$,
then we have an induced map $\overline{\phantom{a}} : V(\lambda)\to V(\lambda)$
by $\overline{u\cdot v_{\lambda}} = \overline{u}\cdot v_{\lambda}$. (This is 
well defined by \cite{jcj}, Proposition 11.9.) The fixed
choice for $v_{\lambda}$ leads to a fixed crystal base $(\LLlam, \B(\lambda))$ 
of $V(\lambda)$, where $\LLlam$ is spanned by all $\Ftilde_{\alpha_{i_1}}\cdots
\Ftilde_{\alpha_{i_r}}(v_{\lambda})$, for $r\geq 0$. Now by, e.g., \cite{lusztig5a}
Theorem 1.8, \cite{jcj} Theorem 11.10, there is a unique basis
$\{ G_{\lambda}(b) \mid b\in \B(\lambda)\}$ of $\LLlam$, such that
\begin{enumerate}
\item $G_{\lambda}(b) = b \bmod q\LLlam$,
\item $\overline{G_{\lambda}(b)}=G_{\lambda}(b)$.
\end{enumerate}
This basis is called the canonical basis of $V(\lambda)$. \par
In the sequel, when we write {\em the} crystal base or {\em the} canonical basis 
of $V(\lambda)$ we always assume that a fixed highest-weight vector $v_{\lambda}$
has been chosen, which makes the choice of crystal base, canonical basis unique.\par
The crystal graph $\Gamma_{\lambda}$ of the module $V(\lambda)$ is defined as follows.
The points of $\Gamma_{\lambda}$ are the elements of $\B(\lambda)$, and there is an
edge $b_1 \stackrel{\alpha}{\longrightarrow} b_2$ if $\Ftilde_{\alpha}(b_1)=b_2$.
There is a very elegant method to compute the crystal graph, using Littelmann's
path method.
Let $\R P$ be the vector space over $\R$ spanned by the weights. Let $\Pi$ be the set
of piecewise linear paths $\pi : [0,1] \to \R P$, such that $\pi(0)=0$.
For $\alpha\in\Delta$ Littelmann
defined operators $e_{\alpha},f_{\alpha}: \Pi \to \Pi \cup \{0\}$ (cf. \cite{litt6},
\cite{litt4}), with the following property. Let $\lambda \in P^+$ be a dominant weight,
and let $\pi_{\lambda}$ be the path given by $\pi_{\lambda}(t)=\lambda t$ (i.e., a 
straight line from the origin to $\lambda$). 
let $\Pi_{\lambda}$ be the set of all $f_{\alpha_{i_1}}\cdots f_{\alpha_{i_k}}(
\pi_{\lambda})$. Then all paths in $\Pi_{\lambda}$ end in an element of $P$. 
Furthermore, the number of
paths ending in $\mu\in P$ is equal to the dimension of the weight space with weight
$\mu$ in the irreducible $U_q$-module $V(\lambda)$. 
Now we consider the directed labeled graph with point set $\Pi_{\lambda}$, and edges
$\pi_1 \stackrel{\alpha}{\longrightarrow} \pi_2$ if $f_{\alpha}(\pi_1)=\pi_2$. This 
graph is isomorphic to the crystal graph of $V(\lambda)$ (\cite{kashiwara2}). \par
Let $M_1,M_2$ be $U_q$-modules, then $M_1\otimes M_2$ is a $U_q$-module via
the comultiplication of $U_q$. There are many possible ways to define this, and
the comultiplication $\Delta : U_q\to U_q\otimes U_q$ that we use is given by
\begin{align*}
\Delta(E_{\alpha}) &= E_{\alpha}\otimes K_{\alpha}^{-1} +1\otimes E_{\alpha}\\
\Delta(F_{\alpha}) &= F_{\alpha}\otimes 1 + K_{\alpha}\otimes F_{\alpha}\\
\Delta(K_{\alpha}) &= K_{\alpha}\otimes K_{\alpha}
\end{align*}
(see \cite{jcj}, 9.13). \par

\section{Canonical bases of tensor products}\label{psi}

Here we give a description of the canonical basis of a tensor product, following
\cite{jcj} Chapter 9, \cite{lusztig} 27.3. \par
Let $V(\mu)$, $V(\mu')$ be two irreducible $U_q$-modules, with highest weights 
$\mu$, $\mu'$. 
Let $C=\{v_1,\ldots ,v_m\}$ and $C'=\{v_1',\ldots ,v_n'\}$ be fixed 
canonical bases of $V(\mu),V(\mu')$. Denote the weight of $v_i$, $v_j'$ by $\nu_i$, 
and $\nu_j'$ 
respectively. Then $\nu_i = \mu -\sum_k a_{k,i} \alpha_k$ with $a_{k,i}\in 
\Z_{\geq 0}$, and we say that $\sum_k a_{k,i}$ is the height of $\nu_i$. The height
of $\nu_j'$ is defined similarly. We assume that the bases $C,C'$ are ordered according
to increasing height. So $v_1 = v_{\mu}$, $v_1'=v_{\mu'}$ are the highest-weight 
vectors. \par
Let $(\LL,\B)$ and $(\LL',\B')$ be crystal bases of $V(\mu)$ and $V(\mu')$ respectively.
Here $\LL$, $\LL'$ are spanned by $C$ and $C'$ respectively. Furthermore, $\B$, $\B'$
consist of the cosets $v_i \bmod q\LL$, $v_j' \bmod q\LL'$. Now by \cite{jcj}, 
Theorem 9.17, $(\LL\otimes \LL',\B\otimes \B')$ is a crystal base of 
$V(\mu)\otimes V(\mu')$. \par
We let $\Theta$ be the element from \cite{lusztig}, 4.1, and $P : U_q\otimes U_q
\to U_q\otimes U_q$ is the algebra homomorphism defined by $P(a\otimes b)=b\otimes a$.
We set $\Theta^0 = P(\overline{\Theta})$; then $\Theta^0 =\sum_{\eta\geq 0} 
\Theta^0_{\eta}$, where the sum runs over all $\eta = \sum_k b_k\alpha_k$ with
$b_k\in \Z_{\geq 0}$. Furthermore, $\Theta^0_{\eta}\in U^+_{\eta}\otimes U^-_{\eta}$,
and $\Theta^0_0 = 1\otimes 1$.
Now $\Psi_0 : V(\mu)\otimes
V(\mu')\to V(\mu)\otimes V(\mu')$ is the map defined by $\Psi_0(v\otimes v') = 
\Theta^0(\overline{v}\otimes \overline{v}')$.

\begin{lem}\label{lem1}
We have $\Psi_0(u\cdot v\otimes v')=\overline{u}\cdot \Psi_0(v\otimes v')$ for all
$u\in U^-$. Furthermore, $\Psi_0^2(v\otimes v') = v\otimes v'$ for all $v\in V(\mu)$,
$v'\in V(\mu')$.
\end{lem}

\begin{pf}
This is the same as the corresponding results in \cite{lusztig}, 27.3.1. The
difference is that we use a different comultiplication. Denoting the comultiplication
used in \cite{lusztig} by$\Delta_L$, we have $\Delta_L(F_{\alpha})= F_{\alpha}\otimes 
K_{\alpha}^{-1}+1\otimes F_{\alpha}$. This means that for $u\in U^-$, we have 
$\Delta(u) = P(\overline{\Delta}_L(u))$, where $\overline{\Delta}_L(u) = 
\overline{\Delta_L(\overline{u})}$. The property $\Delta_L(u)\Theta = \Theta 
\overline{\Delta}_L(u)$ (\cite{lusztig}, Theorem 4.1.2) now translates to 
$\Delta(u)\Theta^0 = \Theta^0 \overline{\Delta}(u)$, where $\overline{\Delta}$
is defined similarly to $\overline{\Delta}_L$. Form this the first statement follows.
The second follows from $\Theta^0\overline{\Theta}^0 = 1\otimes 1$ (\cite{lusztig},
Corollary 4.1.3).  
\end{pf}

We define a partial order on the $v_i\otimes v_j'$. We set 
$v_i\otimes v_j' < v_k\otimes v_l'$ if and only if $i<k$, $j>l$, and $\nu_i+\nu_j'
=\nu_k+\nu_l'$.

\begin{prop}\label{prop1}
There are unique elements $w_{ij}\in V(\mu)\otimes V(\mu')$ such that 
\begin{enumerate}
\item $\Psi_0(w_{ij}) = w_{ij}$,
\item $w_{ij} = v_i\otimes v'_j +\sum_k \zeta_k v_{i_k}\otimes v_{j_k}'$, with
$\zeta_k\in q\Z[q]$, and $\nu_{i_k}+\nu_{j_k}'= \nu_i+\nu'_j$.
\end{enumerate}
Also, $v_{i_k}\otimes v_{j_k}' < v_i\otimes v_j'$ for all $k$.
The elements $w_{ij}$ form a basis of $V(\mu)\otimes V(\mu')$. 
\end{prop}

\begin{pf} This goes in the same way as \cite{lusztig}, Theorem 27.3.2. Note that
\begin{equation}\label{eq1}
\Psi_0(v_i\otimes v_j)= v_i\otimes v_j +\sum_k \xi_k v_{i_k}\otimes v_{j_k}', 
\end{equation}
with $\xi_k\in \Z[q,q^{-1}]$. From $\Theta^0_{\eta}\in U^+_{\eta}\otimes U^-_{\eta}$
and the assumption on the ordering of $C,C'$ it follows that $v_{i_k}\otimes v_{j_k}' 
< v_i\otimes v_j'$ for all $k$. Let $X$ be the set of all $(i,j)$
with $\nu_i+\nu_j' = \nu$, for a certain $\nu$. Order the elements of $X$ in such a way 
that 
$v_i\otimes v_j' < v_k\otimes v_l'$ implies that $(i,j)$ appears before $(k,l)$. 
Let $(i,j)$ be the smallest element of $X$. Then by 
(\ref{eq1}), we see that $\Psi_0(v_i\otimes v_j')= v_i\otimes v_j'$. So in this case
we set $w_{ij}= v_i\otimes v_j'$. Now choose a $(k,l)\in X$, and suppose 
that $w_{r,s}$ exist for all $(r,s)\in X$ appearing before $(k,l)$. Then using 
(\ref{eq1}), and the
triangular form of the $w_{r,s}$ we can write $\Psi_0(v_k\otimes v_l')-v_k\otimes v_l' =
\sum_{r,s} \zeta_{r,s} w_{r,s}$, where $v_r \otimes v_s' < v_k\otimes v_l'$.
After taking images under $\Psi_0$, and using the fact that $\Psi_0$ is
an involution, we see that the $\zeta_{r,s}\in \Z[q,q^{-1}]$ satisfy 
$\bar{\zeta}_{r,s} = -\zeta_{r,s}$. This implies that there are unique $\delta_{r,s}
\in q\Z[q]$ with $\zeta_{r,s}=\delta_{r,s} -\bar{\delta}_{r,s}$. Now set 
$w_{k,l} = v_k\otimes v_l' + \sum_{r,s} \delta_{r,s} w_{r,s}$. For the uniqueness suppose
that there are $w_{ij}'\in V(\mu)\otimes V(\mu')$ satisfying 1., 2. Then we write 
$w_{ij}'$ as a linear combination of $w_{ij}$. By 2. the coefficients are in $\Z[q]$.
Then 1. implies that they are in $\Z$. Finally, from 2. we see that one coefficient is
$1$, and the others are $0$.
\end{pf}

Let $V(\lambda)$ denote the $U_q$-submodule of
$V(\mu)\otimes V(\mu')$ generated by $v_{\mu}\otimes v_{\mu'}=v_1\otimes v_1'$. So 
$V(\lambda)$ is the irreducible $U_q$-module with highest weight $\lambda=\mu+\mu'$.
Set $\LLlam = (\LL\otimes \LL') \cap V(\lambda)$, and $\B(\lambda) = (\B\otimes \B')
\cap \LLlam/q\LLlam$. Then by \cite{jcj}, Proposition 9.10, $(\LLlam, \B(\lambda))$ is
a crystal base of $V(\lambda)$ (the hypotheses of this proposition are satisfied
by \cite{jcj}, Proposition 9.23, Lemma 9.26). 

\begin{thm}\label{thm1}
The elements of the canonical basis of $V(\lambda)$ have the form 
$v_i\otimes v_j' +\sum_k \zeta_k v_{i_k}\otimes v_{j_k}'$, with
$\zeta_k\in q\Z[q]$,  and $v_{i_k}\otimes v_{j_k}' < v_i\otimes v_j'$ for all $k$.
\end{thm}

\begin{pf}
We have that $\Psi_0(v_1\otimes v_1')=v_1\otimes v_1'$. So by Lemma \ref{lem1}, 
$\Psi_0$ coincides with $\overline{\phantom{a}}$ on $V(\lambda)$ (where 
$\overline{u\cdot v_1\otimes v_1'}=\overline{u}\cdot v_1\otimes v_1'$). Hence 
the elements of the 
canonical basis of $V(\lambda)$ are invariant under $\Psi_0$. 
Also, since the elements of the
canonical basis lie in $\LLlam$ and are equal to a $v_i\otimes v_j'\bmod q\LLlam$, 
they must be of the form  $v_i\otimes v_j +\sum_k \zeta_k v_{i_k}\otimes v_{j_k}'$
with $\zeta_k\in q\Z[q]$. 
Now Proposition \ref{prop1} finishes the proof.
\end{pf}

Let $V(\mu_1),\ldots,V(\mu_r)$ be irreducible $U_q$-modules with canonical
bases $C_i =\{v_1^i,\ldots, v_{m_i}^i\}$, ordered according to increasing height.
We consider the tensor product $V=V(\mu_1)\otimes \cdots\otimes V(\mu_r)$. We write
$v_{i_1}^1\otimes \cdots\otimes v_{i_r}^r \Plexs v_{j_1}^1\otimes \cdots\otimes 
v_{j_r}^r$ if there is a $k$ with $i_1=j_1,\ldots,i_k=j_k$ and $i_{k+1} < j_{k+1}$.
Set $\lambda=\mu_1+\cdots +\mu_r$ and let $V(\lambda)$ be the $U_q$-submodule of $V$
generated by $v_1^1\otimes \cdots \otimes v_1^r$. 

\begin{cor}\label{cor2.1}
The elements of the canonical basis of $V(\lambda)$ have the form 
$v_{i_1}^1\otimes \cdots\otimes v_{i_r}^r +\sum_k \zeta_k x_k$
where $\zeta_k\in q\Z[q]$, $x_k\in C_1\otimes \cdots \otimes C_r$, and $x_k \Plexs
v_{i_1}^1\otimes \cdots\otimes v_{i_r}^r$.
\end{cor}

\begin{pf}
The case $r=2$ is covered by Theorem \ref{thm1}, so suppose $r>2$. 
Let $W$ be the $U_q$-submodule of $V(\mu_2)\otimes\cdots\otimes V(\mu_r)$ generated by 
$v_1^2\otimes \cdots \otimes v_1^r$. Then $W$ is the irreducible $U_q$-module with
highest weight $\mu_2+\cdots +\mu_r$. Let $\{w_1,\ldots, w_s\}$ be the canonical 
basis of $W$. Then by Theorem \ref{thm1} the elements of the canonical basis of 
$V(\lambda)$ have the form 
$$v_{i_1}^1\otimes w_{j_1} +\sum_{k\geq 2} \zeta_k v_{i_k}^1\otimes w_{j_k},$$
with $i_k < i_1$ for all $k\geq 2$, and  $\zeta_k\in q\Z[q]$.
We get the result by writing all $w_{j_k}$ for $k\geq 1$ as linear combinations
of elements of $C_2\otimes \cdots \otimes C_r$ and use induction.
\end{pf}

\section{A monomial basis of $V(\lambda)$}\label{monomial}

In this section we first describe a basis of $V(\lambda)$, following \cite{lakshmibai}.
Then using this we derive an algorithm for constructing the canonical basis of $V(\lambda)$,
when $V(\lambda)$ is viewed as a submodule of a tensor product.\par
Let $\pi \in \Pi_{\lambda}$. Then the first direction of $\pi$ is
$w(\lambda)$ for some $w\in W/W_{\lambda}$ (\cite{litt6}, 5.2), where $W_{\lambda}$ is
the stabilizer of $\lambda$. Set $\phi(\pi)=w$. Let $s_{i_1}\cdots s_{i_r}$
be the reduced expression for $\phi(\pi)$, which is lexicographically the smallest.
(Here $s_{i_1}\cdots s_{i_r}$ is lexicographically smaller than $s_{j_1}\cdots s_{j_r}$ 
if there is a $k>0$ such that $i_1=j_1,\ldots, i_{k-1}=j_{k-1}$ and $i_k<j_k$.)
Then we define integers $n_1,\ldots, n_r$,
and paths $\pi_0,\pi_1,\ldots, \pi_r$ in the following way. First, $\pi_0=\pi$.
We let $n_k$ be maximal such that $e_{\alpha_{i_k}}^{n_k}\pi_{k-1} \neq 0$, and we set 
$\pi_k = e_{\alpha_{i_k}}^{n_k}\pi_{k-1}$. Set $\eta_{\pi} = (n_1,\ldots, n_r)$, and
$F_{\pi} = F_{\alpha_{i_1}}^{(n_1)}\cdots F_{\alpha_{i_r}}^{(n_r)}$. Let $b_{\lambda}\in 
\B(\lambda)$ denote the unique element of weight $\lambda$ (it is the coset of
$v_{\lambda}$ modulo $q\LLlam$). Set
$b_{\pi} = \Ftilde_{\alpha_{i_1}}^{n_1}\cdots \Ftilde_{\alpha_{i_r}}^{n_r}(b_{\lambda})$; 
then $\B(\lambda) =\{b_{\pi} \mid \pi \in \Pi_{\lambda} \}$ (this follows from
\cite{kashiwara2}).\par
In the sequel we let $<_B$ denote the Bruhat order on the Weyl group $W$. The 
lexicographical order on sequences of length $r$ is defined by $(m_1,\ldots, m_r)
\Plexs (n_1,\ldots, n_r)$ if there is a $k$ such that $m_1=n_1,\ldots, m_{k-1}=n_{k-1}$ and
$m_k < n_k$. We now define a partial order on $\Pi_{\lambda}$ as follows. First of all, 
$\pi<\sigma$ if  $\phi(\pi)<_B \phi(\sigma)$. Secondly, if 
$\phi(\pi)=\phi(\sigma)$, then $\pi <\sigma$ if $\eta_{\pi} \Plexg \eta_{\sigma}$.
For the proof of the folowing theorem we refer to  \cite{lakshmibai}.

\begin{thm}\label{thm3.1}
$$F_{\pi}\cdot v_{\lambda} = G_{\lambda}(b_{\pi}) + \sum_{\sigma < \pi} 
\zeta_{\pi,\sigma}  G_{\lambda}(b_{\sigma}),$$
where $\zeta_{\pi,\sigma}\in \Z[q,q^{-1}]$.
\end{thm}

\begin{cor}
The set $\{F_{\pi}\cdot v_{\lambda} \mid \pi\in \Pi_{\lambda}\}$ is a basis of 
$V(\lambda)$.
\end{cor}

Let $\pi\in \Pi_{\lambda}$, and $F_{\pi} = F_{\alpha_{i_1}}^{(n_1)}\cdots 
F_{\alpha_{i_r}}^{(n_r)}$. Then we say that $\pi$ is of weight $\nu=\sum_k n_k\alpha_{i_k}$.
We note that this means that $F_{\pi}\cdot v_{\lambda}$ is a weight vector in 
$V(\lambda)$ of weight $\lambda-\nu$. By $\Pi_{\lambda,\nu}$ we denote
the set of all $\pi\in\Pi_{\lambda}$ of weight $\nu$.\par
Suppose that $\lambda = \mu_1+\cdots +\mu_r$, where the $\mu_i$ are dominant weights. 
Also suppose that we are given the modules $V(\mu_i)$ with canonical bases
$C_i= \{v_1^i,\ldots ,v_{m_i}^i\}$, ordered according to increasing height.
We identify $V(\lambda)$ with the $U_q$-submodule of $V(\mu_1)\otimes\cdots \otimes 
V(\mu_r)$ generated by $v_{\lambda} = v_1^1\otimes \cdots \otimes v_r^1$. 
Set $C = C_1\otimes \cdots
\otimes C_r$, which is a basis of $V(\mu_1)\otimes \cdots \otimes V(\mu_r)$ ordered
with respect to $\Plexs$ (see the previous section). \par
Theorem \ref{thm3.1} leads to the following algorithm for computing the 
$G_{\lambda}(b_{\pi})$, for $\pi\in \Pi_{\lambda,\nu}$. Let $\sigma_1,\ldots, \sigma_r$
be the elements from $\Pi_{\lambda,\nu}$ that are smaller than $\pi$. We 
assume that the $G_{\lambda}(b_{\sigma_i})$ already have been computed. Write
$G_{\lambda}(b_{\sigma_i})= y_i + \sum_k \zeta_k y_{i,k}$, where $y_i,y_{i,k}\in C$ and 
$y_i \Plexg y_{i,k}$ for all $k$. We assume that $y_i \Plexs y_j$ implies that $i>j$. 
Then we do the following:
\begin{enumerate}
\item Write $X=F_{\pi}\cdot v_{\lambda}$ as a linear combination of elements
of $C$.
\item For $i=1,\ldots,r$ we do the following. Let $\zeta_i$ be the coefficient of
$y_i$ in $X$. Let $\xi_i$ be the unique element of $\Z[q,q^{-1}]$
such that $\overline{\xi}_i=\xi_i$ and $\zeta_i+\xi_i\in q\Z[q]$. Set $X:= 
X+\xi_iG_{\lambda}(b_{\sigma_i})$.
\end{enumerate}

\begin{prop}
When the loop in Step 2 terminates we have that $X=G_{\lambda}(b_{\pi})$.
\end{prop}

\begin{pf}
Note that by Theorem \ref{thm3.1} there are coefficients $\xi_i$ such that 
$G_{\lambda}(b_{\pi})= F_{\pi} \cdot v_{\lambda} +\sum_{i=1}^r \xi_i G_{\lambda}
(b_{\sigma_i})$. This implies that $\overline{\xi}_i=\xi_i$.
Also, by Corollary \ref{cor2.1}, we have that $G_{\lambda}(b_{\pi})$ is 
of the form $x+\sum_k \omega_k x_k$, where $x,x_k\in C$ and $\omega_k \in q\Z[q]$.
Note that by Corollary \ref{cor2.1} $y_1$ does not occur in any 
$G_{\lambda}(b_{\sigma_i})$, except  $G_{\lambda}(b_{\sigma_1})$. Therefore,
$\xi_1$ is uniquely determined by the requirements that it should be invariant
under $\overline{\phantom{a}}$, and $\zeta_1+\xi_i\in q\Z[q]$. Then in the same way
we see that $\xi_2$ is uniquely determined, and so on.
\end{pf}

\begin{exa}
Let $\Phi$ be the root system of type $G_2$. We denote the simple roots of $\Phi$
by $\alpha$, $\beta$, where $\beta$ is long. The fundamental module $V(\lambda_1)$
is $7$-dimensional, and the canonical basis is $C_1=\{v_1,\ldots, v_7\}$;
these are weight vectors of weights $( 1, 0 )$, $( -1, 1 )$, $( 2, -1 )$, $( 0, 0 )$, 
$( -2, 1 )$, $( 1, -1 )$, $( -1, 0 )$. Here we abbreviate a weight $m\lambda_1+n\lambda_2$
as $(m,n)$. The fundamental module $V(\lambda_2)$ is
$14$-dimensional and has canonical basis $C_2=\{w_1,\ldots, w_{14}\}$. The $w_i$ are
weight vectors of weights $( 0, 1 )$, $( 3, -1 )$, $( 1, 0 )$, $( -1, 1 )$, $( -3, 2 )$, 
$( 2, -1 )$, $( 0, 0 )$, $(0,0)$, $( 3, -2 )$, $( -2, 1 )$, $( 1, -1 )$, $( -1, 0 )$, 
$( -3, 1 )$, $( 0, -1 )$. A description of the action of the generators of $U_q$ on
$V(\lambda_1)$ can for instance be found in \cite{kamis}, and the action of 
$U_q$ on $V(\lambda_2)$ is described in \cite{jcj}, 5A.4. Alternatively, these modules
can be constructed using the {\sf GAP}4 package {\sf QuaGroup} (\cite{gap4}, 
\cite{quagroup}). This package has been used to perform many of the calculations
in the rest of this example. 
Now we set $\lambda=2\lambda_1+\lambda_2$. Then $V(\lambda)$ is
the submodule of $W=V(\lambda_1)\otimes V(\lambda_1)\otimes V(\lambda_2)$ generated
by $v_1\otimes v_1\otimes w_1$. We construct the elements of the canonical basis of 
$V(\lambda)$ that are of weight $\mu=(-2,2)$.  
We use the following elements of weight $\mu$:
\begin{align*}
& x_1 = v_1\otimes v_2\otimes w_{10},~  x_2 = v_1\otimes v_4\otimes w_{5},~ 
x_3 = v_1\otimes v_5\otimes w_{4},~ x_4 = v_2\otimes v_1\otimes w_{10},\\
& x_5 = v_2\otimes v_2\otimes w_{7},~  x_6 = v_2\otimes v_2\otimes w_{8},~ 
x_7 = v_2\otimes v_3\otimes w_{5},~ x_8 = v_2\otimes v_4\otimes w_{4},\\
& x_9 = v_2\otimes v_5\otimes w_{3},~  x_{10} = v_2\otimes v_7\otimes w_{1},~ 
x_{11} = v_3\otimes v_2\otimes w_{5},~ x_{12} = v_4\otimes v_1\otimes w_{5},\\
& x_{13} = v_4\otimes v_2\otimes w_{4},~  x_{14} = v_4\otimes v_5\otimes w_{1},~ 
x_{15} = v_5\otimes v_1\otimes w_{4},~ x_{16} = v_5\otimes v_2\otimes w_{3},\\
& x_{17} = v_5\otimes v_4\otimes w_{1},~  x_{18} = v_7\otimes v_2\otimes w_{1}.
\end{align*} 
They are listed in lexicographical order, i.e., $x_1\Plexs x_2\Plexs \cdots 
\Plexs x_{18}$. The weight space of weight $\mu$ in $V(\lambda)$ is $5$-dimensional.
So we get $5$ paths $\pi_i$ in the crystal graph. The corresponding words in the
Weyl group are $\phi(\pi_1)=s_{\alpha}s_{\beta}s_{\alpha}$, $\phi(\pi_2) = 
s_{\beta}s_{\alpha}s_{\beta}$, $\phi(\pi_3)=s_{\alpha}s_{\beta}s_{\alpha}$,
$\phi(\pi_4)=s_{\alpha}s_{\beta}s_{\alpha}s_{\beta}$, $\phi(\pi_5)=
s_{\alpha}s_{\beta}s_{\alpha}s_{\beta}s_{\alpha}$. Setting $\eta_i=\eta_{\pi_i}$ we have
$\eta_1=(4,2,1)$, $\eta_2=(1,5,1)$, $\eta_3=(3,2,2)$, $\eta_4=(3,1,2,1)$, 
$\eta_5=(2,1,2,1,1)$. So we see that $\pi_1 < \pi_3 < \pi_4 < \pi_5$ and $\pi_2<\pi_4$.
Therefore we have
\begin{align*}
G_{\lambda}(b_{\pi_1}) &= F_{\pi_1}v_{\lambda} = x_{16} + q^2x_{15} + q^3x_{13}+q^6x_{12}+
q^8x_{11}+qx_9 +q^3x_8+q^7x_7+q^5x_3+q^8x_2\\
G_{\lambda}(b_{\pi_2}) &= F_{\pi_2}v_{\lambda} = x_{11}+q^3x_7+q^6x_6.
\end{align*}
Also 
$$F_{\pi_3}v_{\lambda} = x_{17}+q^2x_{16}+q^2x_{14}+q^3x_{13}+q^6x_{11}+q^3x_9+
q^5x_8+q^9x_7.$$
All coefficients, except the first one, are in $q\Z[q]$. Hence $G_{\lambda}(b_{\pi_3}) = 
F_{\pi_3}v_{\lambda}$. Now 
\begin{multline*}
F_{\pi_4}v_{\lambda} = (q+q^{-1})x_{16}+(q+q^3)x_{15}+(1+q^2+q^4)x_{13}+
(q^3+q^5+q^7)x_{12} +  (q^3+q^5+q^7+q^9)x_{11} \\ + (1+q^2)x_9
+(2q^2+q^4)x_8+(q^4+2q^6+q^8)x_7 +q^4x_5+q^6x_4+(q^4+q^6)x_3 +(q^5+q^7+q^9)x_2+q^7x_1.
\end{multline*}
The coefficient of $x_{16}$ is not in $q\Z[q]$. Following the algorithm
we see that $G_{\lambda}(b_{\pi_4}) = F_{\pi_4}v_{\lambda} -(q+q^{-1}) G_{\lambda}(
b_{\pi_1})$; we get
$$G_{\lambda}(b_{\pi_4}) = x_{13}+q^3x_{12}+(q^3+q^5)x_{11}+q^2x_8+(q^4+q^6)x_7+
q^4x_5+q^6x_4+q^5x_2+q^7x_1.$$ 
Finally,
\begin{multline*}
F_{\pi_5}v_{\lambda} = x_{18}+(2q+q^{-1})x_{17}+(2q^3+2q+q^{-1})x_{16}+(q+q^3)x_{15} +
(2q+q^3)x_{14} + (2q^4+3q^2+1)x_{13}\\ +(q^3+q^5+q^7)x_{12}
+(q+2q^3+3q^5+2q^7+q^9)x_{11}+q^3x_{10}+(1+2q^2+2q^4)x_9+(2q^2+3q^4+q^6)x_8 \\
+ (2q^4+3q^6+3q^8+q^{10})x_7 (q^4+q^6)x_5+q^6x_4+(q^4+q^6)x_3+(q^5+q^7+q^9)x_2+q^7x_1.
\end{multline*}
We see that the highest basis vector not having a coefficient in $q\Z[q]$ (apart
from $x_{18}$) is $x_{17}$. So we look at $F_{\pi_5}v_{\lambda} -(q+q^{-1})
G_{\lambda}(b_{\pi_3}) =
x_{18} +(q^3+q+q^{-1})x_{16} +(q^4+2q^2+1)x_{13} +(q^4+q^2+1) + \cdots $
(here all coefficients not written lie in $q\Z[q]$). Now $x_{16}$ does not have a 
coefficient in $q\Z[q]$, so we look at $F_{\pi_5}v_{\lambda} -(q+q^{-1})
G_{\lambda}(b_{\pi_3})
-(q+q^{-1})G_{\lambda}(b_{\pi_1}) = x_{18} + (q^2+1)x_{13} +\cdots $. We see that 
$G_{\lambda}(b_{\pi_5}) = F_{\pi_5}v_{\lambda} -(q+q^{-1})G_{\lambda}(b_{\pi_3})-
(q+q^{-1})G_{\lambda}(b_{\pi_1})-G_{\lambda}(b_{\pi_4})$. 
\end{exa}

\begin{rem}
Let $\pi\in\Pi_{\lambda}$, and let $\phi(\pi)=s_{i_1}\cdots s_{i_r}$ be the 
reduced expression which is the smallest in the lexicographical order. Let
$F_{\pi} = F_{\alpha_{i_1}}^{(n_1)}\cdots F_{\alpha_{i_r}}^{(n_r)}$. Write
$\alpha=\alpha_{i_1}$. If $n_1>1$, then by \cite{litt6}, Lemma 5.3(b), 
$\phi( e_{\alpha}\pi )=\phi(\pi)$, and hence $F_{e_{\alpha}\pi} = 
F_{\alpha_{i_1}}^{(n_1-1)}\cdots F_{\alpha_{i_r}}^{(n_r)}$.
On the other hand, if $n_1=1$, then by \cite{litt6}, Lemma 5.3(a) we see that
$s_{\alpha}\phi(e_{\alpha}\pi)\not<_B \phi(e_{\alpha}\pi)$. So
by \cite{litt6}, Lemma 5.3(b), $\phi(\pi)=s_{\alpha}\phi(e_{\alpha}\pi)$. Therefore
$\phi(e_{\alpha}\pi)=s_{i_2}\cdots s_{i_r}$
which is the smallest (in the lexicographical order) reduced expression for
$\phi(e_{\alpha}\pi)$. Hence $F_{e_{\alpha}\pi} = F_{\alpha_{i_2}}^{(n_2)}\cdots 
F_{\alpha_{i_r}}^{(n_r)}$. The conclusion is that 
$$F_{\pi}\cdot v_{\lambda} = \frac{1}{[n_1]_{\alpha}} F_{\alpha} \cdot ( 
F_{e_{\alpha}\pi} \cdot v_{\lambda}).$$
So in order to compute $F_{\pi}\cdot v_{\lambda}$, we only have to act with 
$F_{\alpha}$ on a vector that we already computed. 
\end{rem}

\begin{rem}
Instead of the algorithm described here for getting the monomials $F_{\pi}$, one can
also follow the procedure outlined in \cite{litt5} for constructing so-called
adapted strings. Instead of $\phi(\pi)$ this procedure uses the longest element in
the Weyl group. However, the monomials one gets in that case are in general different 
from the ones we get. Moreover, in general they do not have the nice property described 
in the previous remark.
\end{rem}

\section{The $A_n$-case}\label{an}
\Yvcentermath1
In this section we assume that the root system $\Phi$ is of type $A_n$.
We use results from \cite{kana} to show that in this case our
algorithm is very much like the algorithm described in \cite{lectof}. \par
The simple roots are $\alpha_1,\ldots, \alpha_n$, where we use the usual ordering
of the nodes of the Dynkin diagram (cf. \cite{bou4}).\par
Since the fundamental weights are all minuscule the corresponding irreducible
$U_q$-modules are easy to construct (cf. \cite{jcj}, Chapter 5A).
For $V(\lambda_k)$ we consider the set of sequences 
$S=\{(i_1,\ldots, i_k)\mid 1\leq i_1 < i_2 <\cdots <i_k\leq n+1\}$. Let $V$ be the 
vectorspace over $\Q(q)$ spanned by $v_s$ for $s\in S$. Let $s\in S$. If $i$ occurs
in $s$, but $i+1$ does not, then let $s^{i-}$ be the sequence obtained from $s$ by
replacing $i$ by $i+1$, and we set $v_s^{i-} = v_{s^{i-}}$. Otherwise $v_s^{i-}=0$. Also, 
if $i+1$ occurs in $s$, but $i$ does not, then let $s^{i+}$ be the sequence obtained from 
$s$ by replacing $i+1$ by $i$, and we set $v_s^{i+} = v_{s^{i+}}$. Otherwise $v_s^{i+}=0$. 
Now a $U_q$-action on $V$ is defined by $F_{\alpha_i} \cdot v_s = v_s^{i-}$, 
$E_{\alpha}\cdot v_s = v_s^{i+}$, and 
$$K_{\alpha_i}\cdot v_s = \begin{cases}
qv_s & \text{ if $i\in S$ and $i+1\not\in S$},\\
q^{-1}v_s & \text{ if $i\not\in S$ and $i+1\in S$},\\
v_s & \text{ otherwise.}
\end{cases}$$
Then the $U_q$-module $V$ is isomorphic to $V(\lambda_k)$. To see this we note that
$v_s$ is a weight vector of weight $\mu_s = a_1\lambda_1+\cdots +a_n\lambda_n$,
where $a_i=1$ if $i\in S$, $i+1\not\in S$, $a_i=-1$ if $i\not\in S$, $i+1\in S$,
and $a_i=0$ otherwise. Set $s_{\lambda_k} = (1,2,\ldots, k)$, then $\mu_{s_{\lambda_k}}=
\lambda_k$. If $i\in S$ and $i+1\not\in S$ then $s_{\alpha_i}(\mu_s) = \mu_{s^{i-}}$.
Since all elements of $S$ can be obtained from $s_{\lambda_k}$ by a sequence of
``moves'' $s \to s^{i-}$ (\cite{kana}, Proposition 3.3.1), we have that $\{\mu_s
\mid s\in S\} = W\cdot \lambda_k$. Finally we compare with \cite{jcj}, 5A.1.\par
Let $\mathcal{L}(\lambda_k)$ be the $A$-submodule of $V$ spanned by the $v_s$, and
let $\B(\lambda_k)$ be the set of all $v_s \bmod q \mathcal{L}(\lambda_k)$. Then
$(\mathcal{L}(\lambda_k), \B(\lambda_k))$ is a crystal base of $V$ 
(\cite{jcj}, Lemma 9.6). Furthermore, $C_k=\{v_s\mid s\in S\}$ is the canonical basis 
of $V$. (Indeed, the $v_s$ are 
certainly invariant under $\overline{\phantom{a}}$ because they are of the form 
$F_{\alpha_{i_1}}\cdots F_{\alpha_{i_r}}\cdot v_{s_{\lambda_k}}$. Secondly $\{v_s\bmod 
q \mathcal{L}(\lambda_k)\mid s\in S\} =\B(\lambda_k)$.)\par
So the elements of $\B(\lambda_k)$ are labeled by the elements of $S$. From \cite{kana}
we get the action of the Kashiwara operators as follows: $\Ftilde_{\alpha_i} (v_s)
=v_s^{i-} \bmod q \mathcal{L}(\lambda_k)$, and $\widetilde{E}_{\alpha_i}(v_s) = 
v_s^{i+}\bmod q \mathcal{L}(\lambda_k)$. \par
We write a sequence $s=(i_1,\ldots, i_k)$ as a diagram with one column of length $k$ 
and the elements
$i_1,\ldots,i_k$ from top to bottom. For example, the sequence $(1,4,5)$ is 
{\tiny $\young(1,4,5)$}. \par
Now let $\lambda=a_1\lambda_1+\cdots + a_n\lambda_n$ be a
dominant weight, and consider the tensor product $W= V(\lambda_1)^{\otimes a_1}
\otimes \cdots \otimes V(\lambda_n)^{\otimes a_n}$. The basis elements are 
labeled by tableaux with $a_n+a_{n-1}+\cdots +a_1$ columns. The first $a_n$ columns
have length $n$, the following $a_{n-1}$ columns have length $n-1$ and so on. The 
tableaux are filled with elements of $\{1,2,\ldots, n+1\}$, such that every column
is strictly increasing. Then every column of length $k$ determines a basis element 
of $V(\lambda_k)$. Tensored together they form a basis element of $W$, e.g.,
$$\young(114,23,3) = \young(4)\otimes \young(1,3)\otimes \young(1,2,3).$$
Then the highest-weight vector $v_{\lambda}$ of weight $\lambda$ in $W$ is labeled by 
the tableau $T_{\lambda}$ with $i$-s in the $i$-th row. Let $V(\lambda)$ 
denote the submodule of
$W$ generated by $v_{\lambda}$. Let $(\LLlam, \B(\lambda))$ be the crystal base
of $V(\lambda)$. Then by \cite{kana}, the elements of $\B(\lambda)$ are labeled 
by tableaux with non-decreasing rows. In particular, these tableaux label the
points in the crystal graph. From \cite{kana} we get 
the following algorithm for computing
$\Ftilde_{\alpha_i}(T)$, $\widetilde{E}_{\alpha_i}(T)$, where $T$ is such a tableau. 
\begin{enumerate}
\item Write the numbers in the tableau as a sequence, starting from the top right,
and going along columns from right to left, top to bottom. Below each number 
write a $+$ if it is equal to $i$, a $-$ if it is $i+1$ and a blank otherwise.
\item If there is a $+$ followed by a $-$ (maybe separated by blanks), then replace
them by blanks. Continue until this operation is no longer possible. 
\item \begin{enumerate} 
\item If there is no $+$ left, then $\Ftilde_{\alpha_i}(T)=0$. Otherwise change the
$i$ corresponding to the leftmost $+$ to a $i+1$. Rebuild the tableau, and the result
is $\Ftilde_{\alpha_i}(T)$.
\item If there is no $-$ left, then $\widetilde{E}_{\alpha_i}(T)=0$. Otherwise
change the $i+1$ corresponding to the rightmost $-$ into $i$. Rebuild the tableau,
the result is $\widetilde{E}_{\alpha_i}(T)$.
\end{enumerate}
\end{enumerate}

\begin{exa}
Let the root system be of type $A_3$, and 
set {\small $T= \young(113,22,3)$}. Then the sequence we get is $3,1,2,1,2,3$. 
If $i=2$, this corresponds to $-o+o+-$ (where we represent a blank by $o$).
After the operation of step 2 this becomes  $-o+ooo$. We see that 
$$\Ftilde_{\alpha_2}(T) = \young(113,23,3),~ \widetilde{E}_{\alpha_i}(T) = 
\young(112,22,3).$$
\end{exa}

The algorithm described in \cite{lectof} for computing the canonical basis of 
$V(\lambda)$ has the same steps as our algorithm: First for every tableau $T$
a monomial $F_T=F_{\alpha_{i_1}}^{(n_1)}\cdots F_{\alpha_{i_t}}^{(r_t)}$ is computed.
Secondly, from the vectors $F_T\cdot v_{\lambda}$ the canonical basis is computed 
using a similar triangular algorithm to the one we use. Therefore, the main difference 
between the algorithms lies in the first step. We investigate this step a little 
further. \par
In \cite{lectof}, 4.1 the authors describe the following algorithm for obtaining a 
monomial $F_T$ from a tableau $T$. Let $i_1$ be the smallest index 
such that $i_1+1$ occurs in an $m$-th row of $T$ with $m\leq i_1$. Furthermore, $r_1$
is the number of occurrences of $i_1+1$ on an $m$-th row with $m\leq i$. Then 
$T_1$ is obtained from $T$ by replacing these $r_1$ occurrences of $i_1+1$ by $i_1$. 
Continuing with $T_1$ instead of $T$ we eventually arrive at the tableau $T_{\lambda}$, 
at which point the algorithm
stops. We have obtained sequences $i_1,\ldots, i_t$, $r_1,\ldots, r_t$ and the monomial
is $F_T=F_{\alpha_{i_1}}^{(n_1)}\cdots F_{\alpha_{i_t}}^{(r_t)}$. \par
Note that applying $\widetilde{E}_{\alpha_i}$ amounts to replacing an $i+1$ by $i$.
Since this $i+1$ was put there by a series of applications of $\Ftilde_{\alpha_{i_k}}$,
starting with $T_{\lambda}$ we see that this $i+1$ must occur on the $m$-th row
with $m\leq i$. By induction on the number
of columns of $T$ it can be shown that if $i_1$ is minimal such that $i_1+1$ occurs in an 
$m$-th row of $T$ with $m\leq i_1$, then $\widetilde{E}_{\alpha_{i_1}}(T)\neq 0$. However
in our algorithm we follow the lexicographically smallest reduced expression of a
word in the Weyl group in order to get the sequence $i_1,\ldots$. This means that 
sometimes we obtain a different monomial than with the algorithm from \cite{lectof},
as the following example shows.

\begin{exa}
Set $T=${\small $\young(114,23,3)$}. Then the monomial obtained by the algorithm of
\cite{lectof} is $F_{\alpha_2}F_{\alpha_3}F_{\alpha_2}F_{\alpha_1}$. Let $\pi$ be
the corresponding path, then $\phi(\pi)=s_{\alpha_3}s_{\alpha_2}s_{\alpha_1}$. This 
means that the monomial that we obtain is $F_{\alpha_3}F_{\alpha_2}^{(2)}F_{\alpha_1}$. 
\end{exa}

We conclude that in the $A_n$-case our algorithm is very similar to, but not the same as,
the algorithm described in \cite{lectof}.

\def\cprime{$'$} \def\cprime{$'$}
\bibliographystyle{plain}

\end{document}